\documentclass[reqno,12pt]{amsart}
\usepackage{amsfonts} %Èç¹ûÔÚ·½À¨ºÅÀïÃæ¼ÓÉÏdraft£¬¿ÉÒÔ¿´µ½ÎÄÖÐÊÇ·ñÓйý³¤µÄ¾ä×Ó¡£
\usepackage{times}
\usepackage{amssymb, amsmath}
\allowdisplaybreaks

\date{July 14, 2009}
\theoremstyle{plain}
\newtheorem{thm}{Theorem}[section]

\theoremstyle{definition}
\newtheorem{defn}{Definition}[section]
\theoremstyle{remark}

\numberwithin{equation}{section}

\title[Axisymmetric Euler-$\alpha$ Equations without Swirl]
{Axisymmetric Euler-$\alpha$ Equations without
Swirl: Existence, Uniqueness, and Radon Measure Valued Solutions}

\author[ Jiu, Niu, Titi and Xin]{}

\begin{document}
\maketitle
%-----------------------------------------------------------------------------------
\centerline{\scshape Quansen Jiu\footnote{The research is partially
supported by National Natural Sciences Foundation of China (No.
10871133) and Project of Beijing Education Committee.}}
\medskip
{\footnotesize
  \centerline{School of Mathematical Sciences, Capital Normal University}
  \centerline{Beijing  100048, P. R. China}
   \centerline{  \it Email: jiuqs@mail.cnu.edu.cn}
   }
\medskip
\centerline{\scshape Dongjuan Niu \footnote{The research is partly supported by
FAPESP grant \# 2008/09473-0.}}
\medskip
{\footnotesize
  \centerline{School of Mathematical Sciences, Capital Normal University}
  \centerline{Beijing  100048, P. R. China}
  \centerline{and}
  \centerline{Departmento de Matematica, IMECC-UNICAMP Caixa Postal 6065,}
  \centerline{Campinas, SP 13081-970, Brasil}
  \centerline{  \it Email: niudongjuan@gmail.com}
}
\medskip
\centerline{\scshape Edriss S. Titi \footnote{The research is supported in part by the
NSF grant no. DMS-0708832 and the ISF grant no. 120/06.}}
\medskip
{\footnotesize \centerline{Department of Computer Science and
Applied Mathematics}
 \centerline{Weizmann Institute of Science, Rehovot 76100,
Israel} \centerline{and}
  \centerline{Department of Mathematics of Mechanical and Aerospace Engineering}
  \centerline{University of California, Irvine, California 92697, USA}
   \centerline{  \it Email: etiti@math.uci.edu}
   }
\medskip
\centerline{\scshape Zhouping Xin\footnote{The research is
partially supported by Zheng Ge Ru Funds, Hong Kong RGC Earmarked
Research Grant CUHK4042/08P and CUHK4040/06P, RGC Central
Allocation Grant CA 05/06. SC01, and a grant from Northwest
University, Xian, China.}}
\medskip
{\footnotesize
  \centerline{The Institute of Mathematical Sciences,
  The Chinese University of Hong Kong}
  \centerline{Shatin, N.T., Hong Kong}
   \centerline{  \it Email: zpxin@ims.cuhk.edu.hk }
   }
%-----------------------------------------------------------------------------------------

\begin{abstract}
The global existence  of weak solutions for the three-dimensional
axisymmetric Euler-$\alpha$  (also known as Lagrangian-averaged
Euler-$\alpha$) equations, without swirl, is established, whenever
the initial unfiltered velocity $v_0$ satisfies $\frac{\nabla \times
v_0}{r}$  is a finite Randon measure with compact support.
Furthermore, the global existence and uniqueness, is also
established in this case provided $\frac{\nabla \times v_0}{r} \in
L^p_c(\mathbb{R}^3)$ with $p>\frac{3}{2}$. It is worth mention that
no such results are known to be available, so far, for the
three-dimensional Euler equations of ideal incompressible flows.


\vskip 0.125in

\end{abstract}

\vskip 0.25in
\begin{center} {\it In honor of Professor John Gibbon on
his 60th birthday}
\end{center}

\section{Introduction}
The Euler-$\alpha$  (also known as the Lagrangian-averaged
Euler-$\alpha$) equations, introduced by Holm, Marsden and Ratiu
\cite{Holm-Marsden-Ratiu}, \cite{HMR-98b},
 are a regularization of the Euler equations
that describe the motion of an ideal incompressible fluid. Moreover,
these same equations also govern the motion of inviscid
(non-viscous) three-dimensional second-grade incompressible
non-Newtonian fluid (see, e.g., \cite{DF1974}, \cite{Dunn95}). The
equations are given by the system:
\begin{align}\label{eqs1}
\left\{
\begin{array}{rl}
&\partial_t v+u\cdot \nabla v+\sum\limits_j v_j \nabla u_j +\nabla p=0,\\
&v=(1-\alpha^2 \Delta) u,\\
&{\rm div}u=0
\end{array}
\right.
\end{align}
in $\mathbb{R}^3$, where $u,p$ represent the velocity vector field and the
 pressure of the fluid, respectively.
%The coefficient $\alpha$ represents the elastic response of the fluid.
The length scale $\alpha>0$ is the regularization parameter, and one recovers,
formally, the Euler equations when $\alpha=0.$

The initial data for \eqref{eqs1} are imposed as
\begin{align}\label{Initial}
u(t=0,x)=u_0(x).
\end{align}

Denoting by $q= {\rm curl} v$, the vorticity, it follows from
\eqref{eqs1} that $q$ satisfies the system:
\begin{align}\label{eqs2}
\left\{
\begin{array}{rl}
&\partial_t q+u\cdot \nabla q=q\cdot \nabla u,\\
&u=K_{\alpha}*q,\\
&q(t=0,x)=q_0(x)
\end{array}
\right.
\end{align}
in $\mathbb{R}^3,$ where $K_{\alpha}$ is the integral kernel of the
inverse of the operator $(1-\alpha^2\Delta) {\rm curl}$. Equation
\eqref{eqs2} resembles the equation of motion of the vorticity in
the three-dimensional Euler equations. However, the vorticity
stretching term, which is the main obstacle for proving the global
regularity for the Euler equations,  is replaced above by the milder
term $q\cdot \nabla u$. Despite this mollification of the vorticity
stretching term the question of global regularity for the
three-dimensional Euler-$\alpha$ is still, as in the case of the
three-dimensional Euler equations, a challenging open problem.

In the following, ${\rm curl} u_0$ and ${\rm curl} v_0$ will denote
the vorticity of $u_0$ and of $v_0=(1-\alpha^2 \Delta) u_0$
respectively. ${\rm curl} v_0$ is also called the initial unfiltered
vorticity, while ${\rm curl}u_0$ is called the filtered initial
vorticity.

There has been a lot of mathematical progress made on the
Euler-$\alpha$ equations recently. For the two-dimensional case,
Kouranbaeva and Oliver \cite{KO} obtained global existence and
uniqueness of \eqref{eqs1}-\eqref{Initial} for the initial
unfiltered vorticity of class $L^2$. The artificial viscosity method
is applied in \cite{KO}. Furthermore, Oliver and Shkoller \cite{OS}
proved the global existence and uniqueness of weak solutions for
initial unfiltered vorticity ${\rm curl}v_0$ in $M(\mathbb{R}^2)$,
which is the space of  finite Radon measures. For the
three-\linebreak dimensional axisymmetric case, Busuioc and Ratiu
\cite{BR} proved the global existence and uniqueness of classical
solutions for the  three-dimensional axisymmetric Euler-$\alpha$
equations without swirl. In particular, in the case of the whole
space $\mathbb{R}^3$, the restrictions on the initial data
 in \cite{BR} are: $u_0\in H^3(\mathbb{R}^3),\
{\rm curl} v_0/r \in L^2(\mathbb{R}^3)$ and ${\rm curl} v_0\in
L^p(\mathbb{R}^3),$ for some $p\in [1,2].$ For the general
three-dimensional case,  a blow up criterion of the smooth
solutions, in the spirit of the Beale-Kato-Majda \cite{BKM}, was
presented in \cite{HL} (see also \cite{Olson-Titi}). The local
existence and a blowup criterion in Besov spaces for the 3D
Lagrangian averaged Euler equations  was given in \cite{LJ}.

Formally, and as we have already mentioned above, when $\alpha=0$ in
\eqref{eqs1}, the Euler-$\alpha$ equations become the classical
Euler equations. A natural question is whether solutions of the
Euler equations can be approximated properly by those of the
corresponding Euler-$\alpha$ equations, especially for the
vortex-sheets initial data. It is well-known that when the  initial
data are a vortex-sheets data, {\it i.e.\ }, the initial vorticity
is a finite Radon measure and the initial velocity is locally
square-integrable, the two-dimensional time-dependent Euler
equations have global (in time) weak solutions when the initial
vorticity ${\rm curl}u_0$ is of one-sign. This result was first
proved by Delort  \cite{D} by regularizing the initial data to
construct the approximate solutions (see also  \cite{LNS} for a
slight generalization). Then the result of Delort was proved by
different approaches (see, e.g.\,, \cite{EM}, \cite{LX},
\cite{Maj2}, \cite{MB}, \cite{Sch}). Specifically, the Navier-Stokes
approximations were used in \cite{Maj2}, \cite{NJX}, and the
vortex-method approximations were applied in \cite{LX}. %Very
%recently, the Euler-$\alpha$ approximations is proved by Bardos,
%Linshiz and Titi (\cite{BLT}). However, the vortex-sheets problem
%for the three-dimensional Euler equations remains unsolved even for
%one-signed case. Along this directions, the readers are referred to
%\cite{Del2} and \cite{JX1}-\cite{JX3} and the further studies will
%be given in \cite{JNXT} by using the Euler-$\alpha$ approximations.
Very recently, the Euler-$\alpha$ equations were  proposed as an
inviscid approximation  to the three-dimensional Euler equations.
Specifically, it was shown in \cite{BLT} and \cite{BLT2} the global
regularity of the vortex sheet problem of the 2D Euler-$\alpha$
equations for a wider class of vortex-sheet. Moreover, it was shown
that these solutions converge to one of the Delort solutions of the
vortex sheet for the 2D Euler equations, as a subsequence of the
regularization parameter
 $\alpha_j\rightarrow 0.$
However, the vortex-sheet problem for the three-dimensional Euler
equations remain unsolved even for the case of one-signed measure.
In a recent work of \cite{BT2} it has been shown, by an example of
3D shear flow of Euler equations, the existence, and persistence for
all time, of singular vortex sheet solutions. This is a
fundamentally different behavior than in the 3D case (see, e.g.,
\cite{CO}, \cite{DR}, \cite{K}, \cite{L}, \cite{W}). For some other
work concerning the vortex sheet problem the reader is referred to
\cite{BT1}, \cite{BC}, \cite{C1}, \cite{C2}, \cite{Del2},
 \cite{JX1}-\cite{JX3}, and further studies
will be reported in the forthcoming paper \cite{JNXT} by using the
Euler-$\alpha$ approximations in the spirit of \cite{BLT} and
\cite{BLT2}.

The purpose of this paper is to prove  the global existence and
uniqueness  of weak solutions for the three-dimensional axisymmetric
Euler-$\alpha$ equations without swirl.   We will first obtain the
global existence and uniqueness of the  solutions when the initial
unfiltered  vorticity $\frac{{\rm curl} v_0}{r}$ belongs to
$L^p_c(\mathbb{R}^3)$ with $p>\frac{3}{2}$, which is the usual $L^p$
Lebesgue space with compact support. Then we will prove the
existence of global weak solutions when the initial unfiltered
vorticity $\frac{{\rm curl} v_0}{r}$ belongs to $M_c(\mathbb{R}^3)$,
which is the space of finite Radon measures with compact support.
The uniqueness of the solutions is still not clear in the case of
weak Radon measures valued solutions. In our analysis,
 the velocity $u(t,x)$ will be recovered from the unfiltered vorticity
$q={\rm curl}\ v$ by the expression  $u=K_{\alpha}*q$, where
$K_{\alpha}$ is the integral kernel of the inverse of
$(1-\alpha^2\Delta) {\rm curl}$ (see \eqref{eqs4} for more details).
As in \cite{BLT}, \cite{BLT2} and \cite{OS}, properties of the
kernel $K_{\alpha}$ near the origin and at the infinity would be
essential to apply the singular integral estimates approach. More
precisely, when $\frac{{\rm curl} v_0}{r}\in L^p_c(\mathbb{R}^3)$
with $p>\frac{3}{2}$, we will be able to prove that the velocity
$u(t,x)$ is Lipschitz continuous with respect to the space variables
and uniformly continuous with respect to the time variable so that
we can prove the global existence and uniqueness of the solution. In
the case where the initial unfiltered vorticity  $\frac{{\rm curl}
v_0}{r}$ belongs to $M_c(\mathbb{R}^3)$ we obtain the global
existence of the weak solutions. This is done  by establishing
appropriate bounds for the gradient and the Hessian of the
approximate solutions and using standard compactness arguments. It
should be noted that in this case the Lipschitz continuity of
$u(t,x)$ on the spatial variables and the uniform continuity with
respect to the time variable are still open.

We will  state our main results in Section 2 and   the proofs of our
main results will be given in Section 3.

\section{Main results}

%Without loss of generality, it is assumed that $\alpha=1$ in this
%paper.
In the cylindrical coordinates $r,\theta,z$, the velocity $u$ is
written as $u=u_re_r+u_{\theta}e_{\theta}+u_ze_z,$ where
$e_r=(\frac{x_1}{r},\frac{x_2}{r},0),
e_{\theta}=(\frac{x_2}{r},-\frac{x_1}{r},0), e_z=(0,0,1), \linebreak
r=(x_1^2+x_2^2)^{\frac{1}{2}}.$ Axisymmetric flows without swirl are
solutions for which the azimuthal (angular) component of the
velocity field satisfies $u_{\theta}\equiv 0$, and $u_r,u_z$ are
independent of $\theta.$ In this case, we also have the unfiltered
velocity $v=v_re_r+v_ze_z$ in the cylindrical coordinate systems. It
should be mentioned that the space of axisymmetric flows is
invariant under the solution of three-dimensional Euler-$\alpha$
equations.

Let $q={\rm curl}\ v$ be the unfiltered vorticity. Then only the
azimuthal (angular) component of $q$ is non-zero in the cylindrical
system, i.e.,
\begin{align}\label{eqs6}
q=q^{\theta}(t,r,z)e_{\theta},
\end{align}
where $q^{\theta}=\partial_z v_r-\partial_r v_z.$
 It follows
that $q$ solves equation \eqref{eqs2} and direct calculations show
that
\begin{align}\label{eqs7}
\partial_t(\frac{q^{\theta}}{r})+u\cdot \nabla
(\frac{q^{\theta}}{r})=0.
\end{align}
That is, the scalar function $q^{\theta}/r$ satisfies the transport
equation. This  is a key and important property of the invariant
family of  axisymmetric flows without swirl (see also \cite{Del2},
\cite{JX1}, \cite{JX2}, \cite{Ettinger-Titi}, \cite{Yudovich}). As
observed above, for this invariant family of flows the vorticity
stretching term, i.e. the right-hand side of equation \eqref{eqs2},
$q\cdot \nabla u \equiv 0$. Consequently, it follows, formally, from
the transport equation \eqref{eqs7} that
\begin{align}
\|\frac{q^{\theta}}{r}\|_{L^r}=\|\frac{q^{\theta}_0}{r}\|_{L^r}\,,
\end{align}
for $r\in [1,\infty]$.

Before stating our main results, we introduce the following
definition of weak solutions. Let $\mathcal{G}$ be the group of all
homeomorphisms $\phi$ of $\mathbb{R}^3$ which preserve the Lebesgue
measure.

\begin{defn}[The case of $L^p(\mathbb{R}^3)(p\ge 1)$]\label{defn1}
Let $T>0$, the vector field $u\in \linebreak
C([0,T];Lip(\mathbb{R}^3))$ and $q={\rm curl}\ v\in
L^p(\mathbb{R}^3)$ with $p\ge 1$ are said to be a weak solution of
\eqref{eqs2} if there exists an unique Lagrangian trajectory
$y(t,x)\in C([0,T]; \mathcal{G})$ satisfying
\begin{align}
&y(t,x_0)=x_0+\int_0^tu(s,y(s,x_0))ds\quad x_0\in \mathbb{R}^3, \label{100}\\
&u(t,x)=K_{\alpha}*q, \label{2-1}\\
&u(t=0,x)=u_0(x), \label{2-2}
\end{align}
and the first equation of $\eqref{eqs2}$ is satisfied in the sense
of distributions. In \eqref{2-1}, $K_{\alpha}$ is, as before,  the
integral kernel representation of the inverse operator of
$(1-\alpha^2\Delta) {\rm curl}$ (see \eqref{eqs4} for more details).
\end{defn}

{\it Remark 2.1} For  the axisymmetric Euler-$\alpha$ equations
without swirl, the first  equation of $\eqref{eqs2}$ becomes
\eqref{eqs7} and in Definition 2.1 it should be satisfied in the
following sense:
\begin{align}
 \int_{\mathbb{R}\times\mathbb{R}^3}(\partial_t
\varphi+u\cdot\nabla \varphi) \frac{q^{\theta}}{r}dxdt=0
\end{align}
 for any $\varphi\in C_0^{\infty}((0,T)\times\mathbb{R}^3)$.

\begin{defn}[The case of $M(\mathbb{R}^3)$]\label{defn2} Let
$T>0$, and let ${\rm curl} v_0/r\in M_c(\mathbb{R}^3)$. The vector
field $u\in L^{\infty}([0,T]\times \mathbb{R}^3)$ is said to be a
weak solution of \eqref{eqs1}   if
\begin{enumerate}
\item
$\nabla u\in L^2([0,T];L^2_{loc}(\mathbb{R}^3)); D^2 u\in
L^1([0,T];L^1_{loc}(\mathbb{R}^3))$.

\item For every test  function $\varphi\in C_0^{\infty}([0,T)\times\mathbb{R}^3)$
with ${\rm div} \varphi=0$, equation \eqref{eqs1} is
 satisfied in the following sense
\begin{align}\label{eqs5}
&\int_{[0,T]\times\mathbb{R}^3}[u(t,x)\cdot(1-\alpha^2\Delta)\partial_t\varphi(t,x)+
(u\cdot\nabla)\varphi\cdot (1-\alpha^2\Delta)u]dxdt\nonumber\\
&+\alpha^2\int_{[0,T]\times\mathbb{R}^3}(\nabla \varphi: D^2)u\cdot
u dxdt =-\int_{\mathbb{R}^3} u_0\cdot
(1-\alpha^2\Delta)\varphi(0,x)dx,
\end{align}
\end{enumerate}
where $\nabla \varphi:
D^2=\sum_{i,k=1}^n\partial_i\varphi_k\partial_k\partial_i$.
\end{defn}

Our main results are stated as:

\begin{thm}\label{thm1}
Assume that the initial velocity is divergence free, axisymmetric
without swirl and ${\rm curl} v_0/r\in L^p_c$ with $p>\frac{3}{2}$.
Then for any $T>0$, there exists a  unique solution  of \eqref{eqs2}
in the sense of  Definition \ref{defn1}, over the interval $[0,T]$.
\end{thm}

\begin{thm}\label{thm2}
Assume that the initial velocity is divergence free, axisymmetric
without swirl and ${\rm curl} v_0/r\in M_c(\mathbb{R}^3).$ Then for
any $T>0$, there exists a global weak solution  $u\in
L^{\infty}([0,T]\times \mathbb{R}^3)$ of \eqref{eqs1} in the sense
of Definition \ref{defn2}. Moreover, we have that $\nabla u\in
L^{\infty}((0,T);L^a+L^\infty)$ with $1\le a<3$ and  $ D^2 u\in
L^{\infty}((0,T);L^b+L^\infty)$
 with $1\le b<\frac{3}{2}$.
\end{thm}

\section{Proof of the Main Results}
In this section, we will give the proofs of Theorems \ref{thm1} and
\ref{thm2}.

For a smooth solution to \eqref{eqs2}, one can define a particle
trajectory  $y(t,x)$ by
\begin{align}\label{eqs25}
&\partial_t y(t,x_0)=u(t,y(t,x_0)),\nonumber\\
&y(0,x_0)=x_0,
\end{align}
where $x_0\in \mathbb{R}^3$. It is noted that the transformation
$y(t,x_0)$ on $\mathbb{R}^3$ preserves the measure due to the
divergence free condition of $u$ (see \cite{MP}). Moreover, one can
recover the velocity $u(t,x)$ from the unfiltered  vorticity $q={\rm
curl}\ v$ through a precise expression of the integral kernel
$K_{\alpha}$ in \eqref{eqs4} as follows. Due to the divergence free
condition of $u,$ there exists a potential vector $\Psi$ such that
${\rm div}\Psi=0$ and $u=\nabla \times \Psi.$ Then by applying the
curl operator to the second equation of $\eqref{eqs1}$ yields that
\begin{equation}\label{3-1}
q=-(1-\alpha^2\Delta)\Delta \Psi.
\end{equation}
Direct calculations  show that the Green function associated with
the operator $(1-\alpha^2\Delta)\Delta$ is (see, e.g., \cite{Kan})
$G_{\alpha}(|x-y|)=\frac{1-e^{-\frac{|x-y|}{\alpha}}}{4\pi |x-y|}.$
Thus, we  deduce that
\begin{align}\label{eqs4}
u(t,x)
&=\nabla\times\int_{\mathbb{R}^3}G_{\alpha}(|x-y|)q(t,y)dy\nonumber\\
&=\int_{\mathbb{R}^3}f_{\alpha}(|x-y|)\frac{x-y}{|x-y|}\times
q(t,y)dy,
\end{align}
where
$f_{\alpha}(|x-y|):=\frac{1}{\alpha^2}f(\frac{|x-y|}{\alpha}),$ and
$f(z):=\frac{(1+z)e^{-z}-1}{4\pi z^2}.$ Obviously, $f(z),\linebreak
\nabla f(z),$ and $zf(z)$ are continuous and bounded functions for
$z\in(0, +\infty)$.
 In addition,  the kernel
$K_{\alpha}$ in the second equation of $\eqref{eqs2}$ can be represented as
 $K_{\alpha}(x,y)=\nabla\times
G_{\alpha}(|x-y|)$.

In view of \eqref{eqs6}, we obtain that for the axisymmetric
Euler-$\alpha$ equations without swirl the velocity is recovered
by
\begin{align}\label{eqs5}
u(t,x)=\int_{\mathbb{R}^3}f_{\alpha}(|x-y|)\frac{x-y}{|x-y|}\times
\frac{q^{\theta}(t,y)}{r}(y_2,-y_1,0)dy.
\end{align}

Now we are ready to prove  Theorem \ref{thm1}. Motivated by
\cite{MP} and \cite{OS}, we will utilize the particle trajectory
\eqref{eqs25} to construct the approximate solutions  and prove the
existence and uniqueness of solutions of the axisymmeric
Euler-$\alpha$ equations without swirl. The key ingredient is to
prove the Lipschitz continuity of the vector fields $u$ with respect
to the spatial variables.

\begin{proof}[Proof of Theorem \ref{thm1}]
The sequences of the approximate solutions of \eqref{eqs7} can be
constructed as follows
\begin{align}
&\partial_t y^n(t,x_0)=u^n(t,y^n(t,x_0)), \label{eqs2-1}\\
&y^n(0,x_0)=x_0, \label{eqs2-2}\\
&y^0(t,x_0)=x_0, \label{eqs2-4}\\
&u^n(t,x)=\int_{\mathbb{R}^3}K_{\alpha}(x,y)q^{n-1}(t,y)dy, \label{eqs2-3}\\
&\frac{(q^n)^{\theta}}{r}(y^n(t,x_0),t)=\frac{q_0^{\theta}}{r}(x_0),\label{eqs2-4+}
\end{align}
for $n\in \mathbb{N}.$

Thanks to \eqref{eqs5}, we have
\begin{align}\label{eqs2-52}
|u^n(t,x)|&= \left|
\int_{\mathbb{R}^3}f_{\alpha}(|x-y^{n-1}(t,z)|)\frac{x-y^{n-1}(t,z)}{|x-y^{n-1}(t,z)|}
\right.\nonumber
\\
&
\times (-x_2+y_2^{n-1}(t,z), x_1-y_1^{n-1}(t,z),0)
\frac{q_0^{\theta}}{r}(z)dz\nonumber\\
&\left. +\int_{\mathbb{R}^3}f_{\alpha}(|x-y^{n-1}(t,z)|)\frac
{x-y^{n-1}(t,z)}{|x-y^{n-1}(t,z)|}\times
(x_2,-x_1,0)\frac{q_0^{\theta}}{r}(z)dz \right|
\nonumber\\
&\le\int_{\mathbb{R}^3}|x-y^{n-1}(t,z)| |f_{\alpha}(|x-y^{n-1}(t,z)|)|
|\frac{q_0^{\theta}}{r}(z)|dz\nonumber\\
&\qquad +\int_{\mathbb{R}^3}|f_{\alpha}(|x-y^{n-1}(t,z)|)|
|\frac{q_0^{\theta}}{r}(z)|dz |x|\nonumber\\
&\leq \frac{C}{\alpha^2}\|\frac{q_0^{\theta}}{r}\|_{L^1}(1+|x|).
\end{align}
In the last inequality above, we have used the facts that
$|x|f_{\alpha}(|x|)$ and $f_{\alpha}(|x|)$ in \eqref{eqs4} are
continuous and bounded functions in $\mathbb{R}^3$. By the
assumption that ${\rm curl} v_0/r\in L^p_c\subset
 L^1$, for $p>3/2$, and \eqref{eqs2-52}, one has
\begin{equation}\label{2-51}
|u^n(t,x)|\le \frac{C}{\alpha^2}(1+|x|),
\end{equation}
where $C$ is a constant depending on the $L^1$ norm of
$\frac{q_0^{\theta}}{r}$ and the bounds of $f_{\alpha}(x)$ and
$xf_{\alpha}(x)$ in \eqref{eqs4}, for $x\in (0,+\infty)$.

Integrating \eqref{eqs2-1} from $0$ to $t$ gives
\begin{align}\label{eqs2-6}
y^n(t,x_0)=x_0+\int_0^tu^n(s,y^n(s,x_0))ds.
\end{align}
\eqref{2-51} and the Gronwall inequality imply that
\begin{align}\label{eqs2-7}
|y^n(t,x_0)|\leq (Ct+|x_0|)e^{Ct}:=L(t,|x_0|),
\end{align}
for any  $t\in [0,T]$. Using \eqref{eqs5}
again yields
\begin{align}\label{eqs2-9}
|u^n(t,x)|&=\left| \int_{\mathbb{R}^3}f_{\alpha}(|x-y^{n-1}(t,z)|)
\frac{x-y^{n-1}(t,z)}{|x-y^{n-1}(t,z)|} \right. \nonumber\\
& \quad \left. \times (y_2^{n-1}(t,z), -y_1^{n-1}(t,z),0)
\frac{q_0^{\theta}}{r}(z)dz \right|\nonumber\\
&\le \frac{C}{\alpha^2}\int_{\mathbb{R}^3}
 L(t,|z|)|\frac{q_0^{\theta}}{r}|(z)dz\nonumber\\
& \le \frac{C}{\alpha^2}  \|\frac{q_0^{\theta}}{r}\|_{L^1},
\end{align}
where  $C$ is a constant depending on $\sup\{|z|:z\in \ {\rm supp}\
\frac{q_0^\theta}{r}\}$.

Next we prove that $u^n$ is Lipshitz continuous with respect to the
spatial variables and is uniformly continuous in time.  For any
$x,x'\in \mathbb{R}^3$, it follows from \eqref{eqs5}-\eqref{eqs2-4+}
that
\begin{align}\label{eqs2-2-1}
&|u^n(t,x)-u^n(t,x^{\prime})|\nonumber\\
&=|\int_{\mathbb{R}^3}\big[f_{\alpha}(|x-y^{n-1}(t,z)|)\frac{x-y^{n-1}(t,z)}
{|x-y^{n-1}(t,z)|}
-f_{\alpha}(|x^{\prime}-y^{n-1}(t,z)|)
\nonumber\\
&\quad \cdot
\frac{x^{\prime}-y^{n-1}(t,z)}{|x^{\prime}-y^{n-1}(t,z)|}\big]
\times (y^{n-1}_2(t,z),-y^{n-1}_1(t,z),0)\frac{q_0^{\theta}}{r}(z)dz|\nonumber\\
&= |\int_{\mathbb{R}^3}
\big[(f_{\alpha}(|x-y^{n-1}(t,z)|)-f_{\alpha}(|x^{\prime}-y^{n-1}(t,z)|))
\frac{x-y^{n-1}(t,z)}{|x-y^{n-1}(t,z)|}
\nonumber\\
&\qquad
+f_{\alpha}(|x^{\prime}-y^{n-1}(t,z)|)\big(\frac{x-y^{n-1}(t,z)}{|x-y^{n-1}(t,z)|}
-\frac{x^{\prime}-y^{n-1}(t,z)}
{|x^{\prime}-y^{n-1}(t,z)|}\big)\big] \nonumber\\
&\qquad \times
(y^{n-1}_2(t,z),-y^{n-1}_1(t,z),0)\frac{q_0^{\theta}}{r}(z)dz\big|.
\end{align}
By the mean value theorem, there exists a point $x''\in
\mathbb{R}^3$ such that
\begin{align}\label{eqs2-2}
&|u^n(t,x)-u^n(t,x^{\prime})|\nonumber\\
 &\leq |x-x^{\prime}|\int_{\mathbb{R}^3} [|\nabla
f_{\alpha}(|x^{\prime\prime}-y^{n-1}(t,z)|)|
+\frac{|f_{\alpha}(|x^{\prime}-y^{n-1}(t,z))|}{|x^{\prime}-y^{n-1}(t,z)|}]\nonumber\\
&\qquad \cdot {|y^{n-1}(t,z)|}
 |\frac{q^{\theta}_0}{r}(z)| dz\nonumber\\
&\leq \frac{c}{\alpha^3}|x-x^{\prime}|\int_{\{z|z\in {\rm
supp}\{\frac{q_0^\theta}{r}\}\}}
(1+\frac{1}{|x^{\prime}-y^{n-1}(t,z)|})
L(t,|z|)\frac{q_0^{\theta}}{r}(z)|dz\nonumber\\
&\leq \frac{c}{\alpha^3}|x-x^{\prime}| (\|\frac{q_0^{\theta}}{r}\|_{L^1}+
\|\frac{q_0^{\theta}}{r}\|_{L^p}\big|\big|\frac{1}{|x^{\prime}-y^{n-1}(t,z)|}
\big|\big|_{L^{\frac{p}{p-1}}})
\nonumber\\
&\leq \frac{c}{\alpha^3}|x-x^{\prime}|
(\|\frac{q_0^{\theta}}{r}\|_{L^1}+\|\frac{q_0^{\theta}}{r}\|_{L^p}),
\end{align}
where $p>\frac 32,$ and $c$ depends on the bounds of  $\nabla
f(|x|)$ and $L(T,|z|)$ (see \eqref{eqs2-7})
if $z\in {\rm supp}\{\frac{q_0^\theta}{r}\}$ and is independent of $\alpha$.

To show the uniform continuity on the time variable, we have, for
any $t, t'\in [0,T]$,
\begin{align}\label{eqs2-10}
&|u^n(t,x)-u^n(t^{\prime},x)|\nonumber\\
&=|\int_{\mathbb{R}^3}\big[f_{\alpha}(|x-y^{n-1}(t,z))|\frac{x-y^{n-1}(t,z)}{|x-y^{n-1}(t,z)|}
\times
(y^{n-1}_2(t,z),-y^{n-1}_1(t,z),0)\nonumber\\
&
-f_{\alpha}(|x-y^{n-1}(t^{\prime},z)|)\frac{x-y^{n-1}(t^{\prime},z)}{|x-y^{n-1}(t^{\prime},z)}
\times (y^{n-1}_2(t^{\prime},z),-y^{n-1}_1(t^{\prime},z),0)\big]
\frac{q_0^{\theta}}{r}(z)dz|\nonumber\\
&\leq \int_{\mathbb{R}^3}\big
|f_{\alpha}(|x-y^{n-1}(t,z)|)-f_{\alpha}(|x-y^{n-1}(t',z)|)|\frac{x-y^{n-1}(t,z)}
{|x-y^{n-1}(t,z)|}|
|y^{n-1}(t,z)| |\frac{q_0^{\theta}}{r}(z)|dz\nonumber\\
&+\int_{\mathbb{R}^3}\big
|f_{\alpha}(|x-y^{n-1}(t',z)|)||\frac{x-y^{n-1}(t,z)}{|x-y^{n-1}(t,z)|}
-\frac{x-y^{n-1}(t',z)}{|x-y^{n-1}(t',z)|}|
|y^{n-1}(t,z)| |\frac{q_0^{\theta}}{r}(z)|
dz\nonumber\\
&+|f_{\alpha}(|x-y^{n-1}(t',z)|)||\frac{x-y^{n-1}(t',z)}{|x-y^{n-1}(t',z)}|
|y^{n-1}(t,z)-y^{n-1}(t',z))||\frac{q_0^{\theta}}{r}(z)|
dz\nonumber\\
 &\leq \|\nabla
f_{\alpha}\|_{L^{\infty}}\int_{\mathbb{R}^3} |y^{n-1}(t,z)-y^{n-1}
(t^{\prime},z)||y^{n-1}(t,z)|
|\frac{q_0^{\theta}}{r}(z)|dz\nonumber\\
 &+\int_{\mathbb{R}^3}|f_{\alpha}(|x-y^{n-1}(t^{\prime},z)|)|
  |\frac{|y^{n-1}(t,z)-y^{n-1}(t^{\prime},z)|}{|x-y^{n-1}
(t^{\prime},z)|} |y^{n-1}(t^{\prime},z)|
|\frac{q_0^{\theta}}{r}(z)|dz\nonumber\\
&+\int_{\mathbb{R}^3}|f_{\alpha}(|x-y^{n-1}(t^{\prime},z)|)|y^{n-1}(t,z)-y^{n-1}(t^{\prime},z)|
|\frac{q_0^{\theta}}{r}(z)|dz\nonumber\\
&\le \frac{C}{\alpha^3}
\int_{\mathbb{R}^3}|y^{n-1}(t,z)-y^{n-1}(t^{\prime},z)|[\frac{1+|x-y^{n-1}
(t^{\prime},z)|}{|x-y^{n-1} (t^{\prime},z)|} L(T,z)+1]
|\frac{q_0^{\theta}}{r}(z)|dz\nonumber\\
&\leq \frac{C}{\alpha^3}|t-t'| \int_{\mathbb{R}^3}|\sup_{s\in
[0,T]}|u^{n-1}(s,z)|[\frac{1+|x-y^{n-1} (t^{\prime},z)|}{|x-y^{n-1}
(t^{\prime},z)|}+1]
|\frac{q_0^{\theta}}{r}(z)|dz\nonumber\\
&\leq
\frac{C}{\alpha^3}|t-t'|(\|\frac{q^{\theta}_0}{r}\|_{L^1}+\|\frac{q_0^{\theta}}{r}\|_{L^p})
\end{align}
for any $p>\frac{3}{2}$, where constant $C$ is independent of $\alpha$ and
 depends on the bounds  of
$u^{n-1}(t,z)$  and $y^{n-1}(t,z)$ if $z$ lies in the support of
$\frac{q_0^\theta}{r}$ in view of \eqref{eqs2-7} and \eqref{eqs2-9}
respectively. Since $u^n$ is Lipschitz continuous with respect to
the spatial variables and is uniformly continuous with respect to
time, and $u^n$ is uniformly bounded by \eqref{eqs2-9}, then the map
$y^n\in C^1([0,T];C(\mathbb{R}^3))$, for every $T>0$. Moreover, we
have that $y^n(t,x)-x\in C^1([0,T];C_B(\mathbb{R}^3))$ in which
$C_B(\mathbb{R}^3)$ means the space of bounded and continuous
functions.

Similar to the estimates in \eqref{eqs2-2} and \eqref{eqs2-10}, one
will prove that $\{y^n(t,x)-x\}$ is a Cauchy sequence in $
C([0,T];C_B(\mathbb{R}^3))$ and furthermore $\{y^n\}$ is a Cauchy
sequence in $C([0,T];\mathcal{G})$. For simplicity, the time
dependence of sequences $\{u^n\}$ and $\{y^n\}, n\in \mathbb{N},$ is
dropped in the following estimates. Note that
\begin{align}\label{eqs2-11}
&|y^n(t,x)-y^{n-1}(t,x)|\nonumber\\
&\leq \int_0^t |u^n(s,y^n(s,x))-u^{n-1}(s,y^{n-1}(s,x))|ds\nonumber\\
&\leq \int_0^t|\int_{\mathbb{R}^3} f_{\alpha}(|y^n(x)-y^{n-1}(z)|)\frac{y^n(x)-y^{n-1}(z)}
{|y^n(x)-y^{n-1}(z)|}\times
(y^{n-1}_2,-y^{n-1}_1,0)\frac{q_0^{\theta}}{r}
\nonumber\\
&\quad
-f_{\alpha}(|y^{n-1}(x)-y^{n-2}(z)|)\frac{y^{n-1}(x)-y^{n-2}(z)}
{|y^{n-1}(x)-y^{n-2}(z)|}\times (y^{n-2}_2,-y^{n-2}_1,0) \frac{q_0^{\theta}}{r}dz|ds
\nonumber\\
&\leq \frac{C}{\alpha^3}\int_0^t \int_{\mathbb{R}^3}
|y^n(x)-y^{n-1}(x)|
|y^{n-1}(z)\frac{q_0^{\theta}}{r}|dz ds\nonumber\\
&\quad +\frac{C}{\alpha^2}\int_0^t\int_{\mathbb{R}^3}
\frac{|y^n(x)-y^{n-1}(x)|}{|y^{n-1}(x)-y^{n-1}(z)|}
 |y^{n-1}| |\frac{q_0^{\theta}}{r}|dzds
\nonumber\\
&\quad
 +
 \frac{C}{\alpha^2}\int_0^t\int_{\mathbb{R}^3}|y^{n-1}-y^{n-2}||\frac{q_0^{\theta}}{r}|
  dzds \nonumber\\
&\leq\frac{C}{\alpha^3}(\|\frac{q_0^{\theta}}{r}\|_{L^1}
+\|\frac{q_0^{\theta}}{r}\|_{L^p})\int_0^t
(\sup\limits_{x\in\mathbb{R}^3}|y^n-y^{n-1}|+
\sup\limits_{x\in\mathbb{R}^3}|y^{n-1}-y^{n-2}|)ds
\end{align}
for $p>\frac 32$, where $C>0$ is a constant depending on
$\sup\{|z|:z\in \ {\rm supp}\ \frac{q_0^\theta}{r}\}$ and is
independent of $\alpha$. Define $g^N(t)=\sup\limits_{n\geq
N}\sup\limits_{x\in\mathbb{R}^3}|y^n(x,t)-y^{n-1}(x,t)|$. Then
\eqref{eqs2-9} implies that $g^N(t)\le \frac{C}{\alpha^3}
||\frac{q_0^\theta}{r}||_{L^1}t<\infty$ for $t\in [0,T]$, where
$C>0$ is a constant depending on $\sup\{|z|:z\in \ {\rm supp}\
\frac{q_0^\theta}{r}\}$. Choose $T_1>0$ small enough such that
$\frac{C}{\alpha^3}(\|\frac{q_0^{\theta}}{r}\|_{L^1}+\|\frac{q_0^{\theta}}{r}\|_{L^p})T_1\le
1/2$. Then it follows  from \eqref{eqs2-11} that
\begin{align}\label{eqs2-12}
g^N(t)\leq k \int_0^t g^{N-1}(s)ds, \ \ t\in [0,T_1],
\end{align}
where $k:=\frac{2
C}{\alpha^3}(\|\frac{q_0^{\theta}}{r}\|_{L^1}+\|\frac{q_0^{\theta}}{r}\|_{L^p})$,
with $C>0$ a constant depending on $\sup\{|z|:z\in \ {\rm supp}\
\frac{q_0^\theta}{r}\}$. Similar as Lemma 3.2 of Chapter 2 in
\cite{MP}, we obtain that
\begin{equation}\label{31}
\lim_{N\rightarrow \infty}g^N(t)\rightarrow 0
\end{equation}
 uniformly on
$[0,T_2]$, for some  $T_2 \in (0, T_1]$ sufficient small. Actually,
it can be proved that $g^N(t)$ can be bounded by the terms of a
convergent geometrical series for $t\in [0,T_2]$. This implies that
 $\{y^n-x\}$ is a Cauchy sequence in
$C([0,T_0];C_B(\mathbb{R}^3))$, where $T_0=\min\{T_1,T_2\}$ depends
on
$\|\frac{q_0^{\theta}}{r}\|_{L^1}+\|\frac{q_0^{\theta}}{r}\|_{L^p}$
and  $\sup\{|z|:z\in \ {\rm supp}\ \frac{q_0^\theta}{r}\}$.  To
continue the procedure above to the interval $[T_0,2T_0]$, we have
\begin{align}
&|y^n(t,x)-y^{n-1}(t,x)|\nonumber\\
&\leq \int_{T_0}^t
|u^n(y^n)-u^{n-1}(y^{n-1})|ds+|y^n(T_0,x)-y^{n-1}(T_0,x)|
\end{align}
for $t\ge T_0$. Due to \eqref{31}, for any $\epsilon>0$, there
exists a $N_1>0$ such that
$\sup\limits_{x\in\mathbb{R}^3}|y^n(T_0,x)-y^{n-1}(T_0,x)|<\epsilon$,
for $n\ge N_1$. Therefore, similar to \eqref{eqs2-11}, for all $n\ge
N_1$, one has
\begin{align}\label{100}
&|y^n(t,x)-y^{n-1}(t,x)|\nonumber\\
&\leq \int_{T_0}^t
|u^n(y^n)-u^{n-1}(y^{n-1})|ds+\epsilon\nonumber\\
&\leq \frac{C}{\alpha^3}(\|\frac{q_0^{\theta}}{r}\|_{L^1}
+\|\frac{q_0^{\theta}}{r}\|_{L^p})\int_{T_0}^t
(\sup\limits_{x\in\mathbb{R}^3}|y^n-y^{n-1}|+\sup\limits_{x\in\mathbb{R}^3}|y^{n-1}-y^{n-2}|)ds+\epsilon,
\end{align}
for $t\in [T_0,2T_0]$. It follows that
\begin{align}
&\int_{T_0}^t\sup\limits_{x\in\mathbb{R}^3}|y^n(s,x)-y^{n-1}(s,x)| ds\nonumber\\
&\leq \frac{C}{\alpha^3}(\|\frac{q_0^{\theta}}{r}\|_{L^1}
+\|\frac{q_0^{\theta}}{r}\|_{L^p})(t-T_0)\int_{T_0}^t
(\sup\limits_{x\in\mathbb{R}^3}|y^n-y^{n-1}|+\sup\limits_{x\in\mathbb{R}^3}|y^{n-1}-y^{n-2}|)ds+\epsilon(t-T_0),\nonumber
\end{align}
for $t\in [T_0,2T_0]$. This implies that
\begin{align}\label{101}
&\int_{T_0}^t\sup\limits_{x\in\mathbb{R}^3}|y^n(s,x)-y^{n-1}(s,x)| ds\nonumber\\
&\leq \int_{T_0}^t
\sup\limits_{x\in\mathbb{R}^3}|y^{n-1}-y^{n-2}|ds+2\epsilon(t-T_0),
\end{align}
Putting \eqref{101} into  \eqref{100}, we get
\begin{align}
&\sup\limits_{x\in\mathbb{R}^3}|y^n(t,x)-y^{n-1}(t,x)|\nonumber\\
&\leq 2\frac{C}{\alpha^3}(\|\frac{q_0^{\theta}}{r}\|_{L^1}
+\|\frac{q_0^{\theta}}{r}\|_{L^p})\int_{T_0}^t
(\sup\limits_{x\in\mathbb{R}^3}|y^n-y^{n-1}|+\sup\limits_{x\in\mathbb{R}^3}|y^{n-1}-y^{n-2}|)ds\nonumber\\
&+2\epsilon(t-T_0)\frac{C}{\alpha^3}(\|\frac{q_0^{\theta}}{r}\|_{L^1}
+\|\frac{q_0^{\theta}}{r}\|_{L^p})+\epsilon.
\end{align}
Since $\frac{C}{\alpha^3}(\|\frac{q_0^{\theta}}{r}\|_{L^1}
+\|\frac{q_0^{\theta}}{r}\|_{L^p})(t-T_0)\le 1/2$, for all $t\in
[T_0,2T_0]$, it follows that
\begin{align}\label{eqs2-12+}
g^N(t)\leq k \int_{T_0}^t g^{N-1}(s)ds+2\epsilon, \ \ t\in [T_0,
2T_0],
\end{align}
where $k$ is same as in \eqref{eqs2-12}. By the arbitrariness of
$\epsilon$, similar to the previous procedure (see also Lemma 3.2 of
Chapter 2 in \cite{MP}), we can prove that $\{y^n-x\}$ is a Cauchy
sequence in $C([T_0,2T_0];C_B(\mathbb{R}^3))$, and furthermore  in
$C([0,T];C_B(\mathbb{R}^3))$ ; and $\{y^n\}$ is a Cauchy sequence in
$C([0,T];\mathcal{G})$, for every $T>0$. The limit of $\{y^n(t,x)\}$
in $C([0,T];\mathcal{G})$ is denoted by $y(t,x)$.

Define $q(t,x)=\frac{q^{\theta}}{r}(x,t)(x_2,-x_1,0)$ satisfying
$\frac{q^{\theta}}{r}(y(t,x),t)=\frac{q_0^{\theta}}{r}(x,t)$. Then
it is straightforward  to prove that
$$q^n(t,x)=\frac{(q^n)^{\theta}}{r}(x,t)(x_2,-x_1,0) \rightharpoonup  q(t,x)$$
in the sense of weakly-* convergence in
$L^\infty([0,T];L^p(\mathbb{R}^3))$, with $p>3/2$, and
$$u^n(t,x)\rightarrow u(t,x)$$
in $C_B([0,T]\times \mathbb{R}^3)$. Moreover, $y, u, q$ solve
\eqref{100}, \eqref{2-1} and satisfy
\begin{align}
 \int_{\mathbb{R}\times\mathbb{R}^3}(\partial_t
\varphi+u\cdot\nabla \varphi) \frac{q^{\theta}}{r}dxdt=0
\end{align}
 for any $\varphi\in C_0^{\infty}((0,T)\times\mathbb{R}^3)$.

The uniqueness can be shown by the  direct estimate on the
difference of two flow maps.  The estimates are similar to the
previous ones and we omit the details here. The proof of the
theorem is finished.
\end{proof}

\begin{proof}[Proof of Theorem \ref{thm2}] Mollifying the initial potential
vorticity,  we can construct the approximate solutions of
\eqref{eqs2} by solving the following problem:
\begin{align}\label{eqs8}
\left\{
\begin{array}{rl}
&\partial_t q^{\epsilon}+u^{\epsilon}\cdot \nabla q^{\epsilon}=q^{\epsilon}\cdot \nabla
u^{\epsilon},\\
&u^{\epsilon}=K_{\alpha}*q^{\epsilon},\\
&q^{\epsilon}(t=0,x)=q_0^{\epsilon},
\end{array}
\right.
\end{align}
where $q_0^{\epsilon}$ is a smooth vector with compact support
which converges to $q_0$ in $M(\mathbb{R}^3)$ and
$\|q_0^{\epsilon}\|_{L^1}\le \|q_0\|_{M}$. Then there exits an
unique smooth solution $(u^\epsilon, q^\epsilon)$ to
\eqref{eqs8} (see \cite{BR} and references therein). Moreover,
there is a smooth pressure $p^{\epsilon}$ such that
$u^\epsilon$ satisfies
\begin{align}\label{eqs101}
\left\{
\begin{array}{rl}
&\partial_t v^\epsilon+u^\epsilon\cdot \nabla v^\epsilon+\sum\limits_j v^\epsilon_j
 \nabla u^\epsilon_j +\nabla p^\epsilon=0,\\
&v^\epsilon=(1-\alpha^2\Delta) u^\epsilon,\\
&{\rm div}u^\epsilon=0.
\end{array}
\right.
\end{align}
 Therefore, for any test
function $\varphi\in C_0^{\infty}([0,T),\mathbb{R}^3)$, satisfying
${\rm div} \varphi=0$, integration by parts yields
 \begin{align}\label{eqs102}
&\int_{[0,T]\times\mathbb{R}^3}[u^{\epsilon}(t,x)(1-\alpha^2\Delta)\partial_t\varphi(t,x)+
(u^{\epsilon}\cdot\nabla)\varphi\cdot (1-\alpha^2\Delta)u^\epsilon]dxdt\nonumber\\
&+\alpha^2\int_{[0,T]\times\mathbb{R}^3}(\nabla \varphi:
D^2)u^\epsilon\cdot u^\epsilon dxdt =-\int_{\mathbb{R}^3}
u^\epsilon_0 (1-\alpha^2\Delta)\varphi(0,x)dx.
\end{align}

Similar analysis as for \eqref{eqs2-9} shows that
\begin{align}\label{eqs26}
\|u^{\epsilon}\|_{L^{\infty}([0,T]\times \mathbb{R}^3)}\leq \frac{C}{\alpha^2}
\|\frac{q_0^{\theta}}{r}\|_{M},
\end{align}
where $C>0$ is a constant depending on $T$ and the Lebesgue measure
of the support of $\frac{q_0^{\theta}}{r}.$
%$\|\frac{q_0^{\theta}}{r}\|_{M}$ and the support of
%$\frac{q_0^{\theta}}{r}$.

We now  estimate  $\nabla u^\epsilon$ and $D^2 u^\epsilon$. Let
$\chi: \mathbb{R} \to [0,1]$ be a smooth functions satisfying
$$
\chi(s)=\left\{
\begin{array}{ll}
&1,\ \ |s|<1,\\[3mm]
&0,\ \ |s|>2 \,.
\end{array}
\right.
$$
 Then it follows from
(3.4) that
\begin{align}
 &\partial_{x_i}
 u^{\epsilon}(t,x)\nonumber\\
 &=\int_{\mathbb{R}^3}\partial_{x_i}[f_{\alpha}(|x-y^{\epsilon}|)
 \frac{x-y^{\epsilon}}{|x-y^{\epsilon}|}]
 (1-\chi(|x-y^{\epsilon}|))\times\frac{q^{\theta}(t,y^{\epsilon})}{r}
 (y_2^{\epsilon},-y_1^{\epsilon},0)dy\nonumber\\
&\quad +\int_{\mathbb{R}^3}\partial_{x_i}[f_{\alpha}(|x-y^{\epsilon}|)
\frac{x-y^{\epsilon}}{|x-y^{\epsilon}|}]
\chi(|x-y^{\epsilon}|)\times
\frac{q^{\theta}(t,y^{\epsilon})}{r}(y_2^{\epsilon},-y_1^{\epsilon},0)dy\nonumber\\
&\equiv G_1(t,x)+G_2(t,x),
\end{align}
for $i=1, 2, 3.$ It is clear that
\begin{align}
|G_1(t,x)|&\le
\frac{C}{\alpha^3}\int_{\mathbb{R}^3}\frac{1}{|x-y^{\epsilon}(t,z)|}
(1-\chi(|x-y^{\epsilon}(t,z)|))
|\frac{(q_0^{\epsilon})^{\theta}}{r}(z)||y^{\epsilon}(t,z)|dz\nonumber\\
&\leq \frac{C}{\alpha^3}\|\frac{q_0^{\theta}}{r}\|_{M},\nonumber
\end{align}
for $(t,x)\in [0,T]\times \mathbb{R}^3.$ One can use Young's
inequality for convolutions to obtain
\begin{align}\label{eqs21-}
&\|G_2(t,x)(t,x)\|_{L^{\infty}((0,T);L^a)} \nonumber\\
&\leq
\sup_{t\in[0,T]}\|\int_{\mathbb{R}^3}\frac{1}{|x-y^{\epsilon}(t,z)|}
\chi(|x-y^{\epsilon}(t,z)|)|\frac{(q_0^{\epsilon})^{\theta}}{r}(z)||y^{\epsilon}(t,z)|dz\|_{L^a}
\nonumber\\
&\leq \frac{C}{\alpha^3}\|\frac{(q_0^{\epsilon})^{\theta}}{r}\|_{L^1}\leq
\frac{C}{\alpha^3}\|\frac{q_0^{\theta}}{r}\|_{M}
\end{align}
for  $1\le a<3$. Thus $\nabla u^\epsilon(t,x)$ is bounded in
$L^{\infty}((0,T);L^a+L^\infty)$ and
\begin{equation}\label{eqs21}
\|\nabla u^{\epsilon}(t,x)\|_{L^{\infty}((0,T);L^a+L^\infty)} \leq
\frac{C}{\alpha^3}\|\frac{q_0^{\theta}}{r}\|_{M}
\end{equation}
for  $1\le a<3$.  Similarly, we can prove that $D^2 u^\epsilon(t,x)$
is bounded in $L^{\infty}((0,T);L^b+L^\infty)$ and that
\begin{equation}\label{eqs21+}
\|D^2 u^{\epsilon}(t,x)\|_{L^{\infty}((0,T);L^b+L^\infty)} \leq
\frac{C}{\alpha^4}\|\frac{q_0^{\theta}}{r}\|_{M}
\end{equation}
for $1\le b< \frac{3}{2}$. In \eqref{eqs21} and \eqref{eqs21+}, the
constant  $C$ is a positive constant depending on
$\|\frac{q_0^{\theta}}{r}\|_{M}$ and the support of
$\frac{q_0^{\theta}}{r}$. In view of \eqref{eqs26},\eqref{eqs21} and
\eqref{eqs21+}, it is easy to obtain   that the terms
$u^\epsilon\cdot \nabla v^\epsilon$, $\sum\limits_j v^\epsilon_j
\nabla u^\epsilon_j$ and $\nabla p^\epsilon$  are bounded in
$L^\infty((0,T);W^{-2,2}_{loc}(\mathbb{R}^3))$ and hence it follows
from \eqref{eqs101} that $\partial_t u^\epsilon$ is bounded in
$L^\infty((0,T);L^2_{loc}(\mathbb{R}^3))$. Note that $u^\epsilon$ is
bounded in $L^\infty((0,T);W^{1,a}_{loc}(\mathbb{R}^3))$ for any
$1\le a<3$. By the Aubin-Lions Lemma (see, e.g.,
\cite{Constantin-Foias}, \cite{Temam}), we obtain that (up to a
subsequence) that $u^\epsilon \rightarrow u$ strongly in
$C([0,T];L^c_{loc}(\mathbb{R}^3)$, with $1\le c<\frac{3a}{3-a}$,
where $u\in L^{\infty}([0,T]\times \mathbb{R}^3)$ is a function
satisfying
 $\nabla u\in L^{\infty}((0,T);L^a+L^\infty)$, with $1\le a<3$, and
   $ D^2 u\in L^{\infty}((0,T);L^b+L^\infty)$,
 with $1\le b<\frac{3}{2}$. Thanks to $\eqref{eqs26}$, we have that
  $u^\epsilon \rightarrow u$ strongly in
$L^\sigma([0,T];L^\sigma_{loc}(\mathbb{R}^3))$, for any $\sigma \in
(1,\infty)$. Moreover, it is clear that $D^2
u^{\epsilon}\rightharpoonup D^2 u$ weakly-* convergence in
$L^{\infty}((0,T);L^b_{loc})$, for any $1\le b<3/2$. Thus, it
follows from \eqref{eqs102} that
\begin{align}
&\int_{[0,T]\times\mathbb{R}^3}[u(t,x)(1-\alpha^2\Delta)\partial_t\varphi(t,x)+
(u\cdot\nabla)\varphi\cdot (1-\alpha^2\Delta)u]dxdt\nonumber\\
&+\alpha^2\int_{[0,T]\times\mathbb{R}^3}(\nabla \varphi:D^2)u\cdot
udxdt =-\int_{\mathbb{R}^3} u_0
(1-\alpha^2\Delta)\varphi(0,x)dx\nonumber
\end{align}
for any  $\varphi\in C_0^{\infty}([0,T),\mathbb{R}^3)$ satisfying
${\rm div} \varphi=0$. The proof of Theorem 2.2 is complete.
\end{proof}


\begin{thebibliography}{100}

\bibitem{BLT}
C. Bardos, J. Linshiz and E.S. Titi.
\newblock Global regularity for a Birkhoff-Rott-$\alpha$ approximation of the
dynamics of vortex sheets of the 2D Euler equations.
\newblock {\em Physica D}, 237:1905--1911, 2008.

\bibitem{BLT2}
C. Bardos, J. S. Linshiz, E. S. Titi.
\newblock Global regularity and convergence of a Birkhoff-Rott-$\alpha$
approximation of the dynamics of vortex sheets of the 2D Euler
equations.
\newblock {\em Comm. Pure Appl. Math.}, (to appear).

\bibitem{BT1}
C. Bardos and E.S. Titi.
\newblock Euler equations of
incompressible ideal fluids.
\newblock {\em Uspekhi Matematicheskikh Nauk, UMN,}
62:3(375):5-6, 2007.  Also in {\em Russian Mathematical Surveys,}
62(3):409--451, 2007.

\bibitem{BT2}
C. Bardos and E.S. Titi.
\newblock Loss of smoothness and energy conserving rough weak solutions for the $3d$
Euler equations.
\newblock {\em Discrete and Continuous Dynamical Systems}, (submitted), 2009.

\bibitem{BKM}
J.T. Beale, T. Kato and A. Majda.
\newblock Remarks on the breakdown of smooth solutions for the 3D Euler equations.
\newblock{\em Commun. Math. Phys.}, 94: 61--66, 1984.

\bibitem{BC}
A. L. Bertozzi, P. Constantin.
\newblock Global regularity for vortex patches.
\newblock{\em Comm. Math. Phys.}, 152(1):19--28, 1993.

\bibitem{BR}
A. V. Busuioc, T. S. Ratiu.
\newblock Some remarks on a certain class of axisymmetric fluids of differential type.
\newblock {\em Physical D}, 191:106--120, 2004.

\bibitem{CO}
R. Caflisch and O. Orellana.
\newblock Singular solutions  and ill-posedness for the evolution of vortex sheets.
\newblock {\em  SIAM J. Math. Anal.},  20(2):293--307, 1989.


\bibitem{C1}
J.-Y. Chemin.
\newblock Persitance de Structures geometriques dans les fluids
incompressibles bidimensionnels.
\newblock{\em Ann. Sci. Ecole Norm. Sup.}, 26(4):517--542, 1993.

\bibitem{C2}
J. Y. Chemin.
\newblock Two-dimensional Euler system and the vortex patches problem.
\newblock {\em In Handbook of mathematical fluid dynamics,}
III:83-160, North-Holland, 2004.

\bibitem{Constantin-Foias} P. Constantin and C. Foias.
\newblock The Navier-Stokes Equations.
\newblock{\em The University of Chicago Press}, Chicago, 1988.

\bibitem{D}
J.-M. Delort.
\newblock Existence de nappes de tourbillon en dimension deux.
\newblock{\em J. Amer. Math. Soc.},4(3):553--586, 1991.

\bibitem{Del2}
J. M. Delort.
\newblock Une remarque sur le probleme des nappes de tourbillon
axisymetriques sur $R^3$.
\newblock {\em J. Funct. Anal.}, 108:274--295, 1992.

\bibitem {DR}
J. Duchon and R. Robert.
\newblock Global vortex sheet solutions of
Euler equations in the plane.
\newblock {\em J. Differential Equations}, 73(2):215--224, 1988.

\bibitem{DF1974}
J. E. Dunn and R. L. Fosdick.
\newblock Thermodynamics, stability, and
boundedness of fluids of complexity 2 and fluids of second grade.
\newblock {\em   Arch. Rat. Mech. Anal.}, 56:191--252, 1974.

\bibitem{Dunn95}
J. E. Dunn and K.R. Rajagopal.
\newblock Fluids of differential type:
    critical reviews and thermodynamic analysis.
\newblock {\em  Int. J. Engng. Sci.},  33:689--729, 1995.

\bibitem{Ettinger-Titi} B. Ettinger and E.S. Titi.
 \newblock Global existence and uniqueness
of weak solutions of 3-D Euler equations with helical symmetry in
the absence of vorticity stretching.
\newblock {\em SIAM J. Math. Anal.},  41(1):269–-296, 2009.


\bibitem{EM}
L. C. Evans and S. M\"uller.
\newblock Hardy space and the two-dimensional Euler equations with non-negative vorticity.
\newblock {\em J. Amer. Math. Soc.}, 7:199--219, 1994.







\bibitem{Holm-Marsden-Ratiu}  D. D. Holm, J. E. Marsden and T. Ratiu.
\newblock Euler-Poincar\'{e} models of ideal fluids with nonlinear
dispersion.
\newblock {\em Phys. Rev. Lett.},  80:4173--4176, 1998.

\bibitem{HMR-98b}   D. D. Holm, J. E. Marsden, and T.S. Ratiu.
\newblock Euler-Poincar\'e equations and semidirect products with applications
to continuum theories.
\newblock {\em Adv. in Math.}, 137:1--81, 1998.

\bibitem{HL}
T. Y. Hou, C. Li.
\newblock On global well-posedness of the lagrangian averaged Euler equations.
\newblock {\em SIAM J. Math. Anal.}, 38(3):782--794, 2006.

\bibitem{JNXT}
Q. S. Jiu, D. J. Niu, E. S. Titi, Z. P. Xin.
\newblock The Euler-$\alpha$
approximations to the 3D axisymmetric Euler equations with
vortex-shhets initial data, preprint, 2009.



\bibitem{JX1}
Q. S. Jiu, Z. P. Xin.
\newblock Viscous approximation and decay rate of  maximal vorticity function
for 3-D axisymmetric Euler equations.
\newblock {\em Acta Math. Sin. (Engl. Ser.)}, 20(3):385--404, 2004.


\bibitem{JX2}
Q. S. Jiu and Z. P. Xin.
\newblock On strong convergence to 3-{D} axisymmetric vortex sheets.
\newblock {\em J. Differential Equations}, 223(1):33--50, 2006.

\bibitem{JX3}
Q. S. Jiu and Z. P. Xin.
\newblock  Smooth Approximations and Exact Solutions of the 3D Steady
Axisymmetric Euler Equations.
\newblock {\em Comm. Math. Phys.}, 287:323--349, 2009.


\bibitem{Kan}
R. P. Kanwal.
\newblock Generalized Functions Theory and Technique.
\newblock {\em Academic Press,} 1983.

\bibitem{KO}
S. Kouranbaeva, M. Oliver.
\newblock Global well-posedness for the averaged Euler equations in two dimensions.
\newblock {\em Physica D}, 138:197--209, 2000.

\bibitem{K}
 R. Krasny.
\newblock Computation of vortex sheet roll-up in
  Trefftz plane.
\newblock {\em J. Fluid Mech.}, 184:123--155, 1987.

\bibitem{L}
G. Lebeau.
\newblock R\'egularit\'e du probl\`eme de Kelvin-Helmholtz
pour l'\'equation d'Euler 2d.
\newblock {\em ESAIM Control Optim. Calc. Var.}, 8:801--825, 2002 (electronic).



\bibitem{LJ}
X. F. Liu and H. Y. Jia.
\newblock Local existence and blowup criterion of the Langrangian
averaged Euler equations in Besov spaces.
\newblock{\em Commun. Pure Appl. Anal.}, 7(4):845--852, 2008.


\bibitem{LX}
J. G. Liu and Z. P. Xin.
\newblock Convergence of vortex methods for weak
solutions to the 2-D Euler equations with vortex sheet data.
\newblock{\em Comm. Pure Appl. Math.}, 48:611--628, 1995.



\bibitem{LNS}
M. C. Lopes Filho, H. J. Nussenzveig Lopes, S. Schochet.
\newblock A criterion for the equivalence of the Birkhoff--Rott and Euler
description of vortex sheet evolution.
\newblock{\em  Trans. Amer. Math. Soc.},359:4125--4142, 2007 (electronic).

\bibitem{Maj2}
A. Majda.
\newblock Remarks on weak solutions for vortex sheets with a distinguished sigh.
\newblock {\em Indiana Univ. Math. J.}, 42:921--939, 1993.

\bibitem{MB}
A. Majda and A. Bertozzi.
\newblock  Vorticity and incompressible flow. Volume~27 of
Cambridge  Texts in Applied Mathematics.
\newblock {\em Cambridge University Press}, Cambridge, 2002.


\bibitem{MP}
C. Marchioro, M. Pulvirenti.
\newblock Mathematical Theory of Incompressible Nonviscous Fluids.
\newblock{\em Springer-Verlag,} New York, 1994.

\bibitem{NJX}
D. J. Niu, Q. S. Jiu and Z. P. Xin.
\newblock Navier-Stokes approximations to 2D vortex sheets in half plane.
\newblock {\em Methods and Applications of Analysis },14(3):263--272, 2007.

\bibitem{OS}
M. Oliver, S. Shkoller.
\newblock The vortex blob method as a second-grade non-newtonian fluid.
\newblock {\em Commun. in Partial Differential Equations}, 26(1\&2):295--314, 2001.

\bibitem{Olson-Titi}
E. Olson and E.S. Titi.
\newblock Viscosity versus
vorticity stretching: global well-posedness for a family of the
Navier-Stokes alpha-like models.
\newblock {\em Nonlinear Analysis}, 66(11):2427--2458, 2007.

\bibitem{Sch} S. Schochet.
\newblock The weak vorticity formulation of the 2D Euler equations
and concentration-cancellation.
\newblock {\em Comm. P. D. E.}, 20:1077--1104, 1995.

\bibitem{Temam} R. Temam.
\newblock Navier-Stokes Equations, Theory and Numerical Analysis. 3rd
revised edition,
\newblock{\em  North-Holland}, 2001.


\bibitem{SY}
T. Shirota, T. Yanagisawa.
\newblock Note on global existence for axially symmetric solutions of the Euler system.
\newblock {\em Proc. Japan. Acad.}, 70(A):299--304, 1994.

\bibitem{Yudovich}
V.I. Yudovich.
\newblock Non-stationary flows of an ideal incompressible fluid.
\newblock{\em Z. Vy\v cisl. Mat. i Mat. Fiz.}, 3:1032--1066, 1963.


\bibitem{W}
S. Wu.
\newblock  Mathematical analysis of vortex sheets.
\newblock {\em Comm. Pure Appl. Math.}, 59(8):1065--1206, 2006.


\end{thebibliography}
\end{document}